\documentclass[amstex,12pt,russian,amssymb]{article}

\usepackage{mathtext}
\usepackage[cp1251]{inputenc}
\usepackage[T2A]{fontenc}
\usepackage[russian]{babel}
\usepackage[dvips]{graphicx}
\usepackage{amsmath}
\usepackage{amssymb}
\usepackage{amsxtra}
\usepackage{latexsym}
\usepackage{ifthen}

\begin{document}

\newcounter{lemma}
\newcommand{\lemma}{\par \refstepcounter{lemma}%
{\bf Лема \arabic{lemma}.}}

\newcounter{corollary}
\newcommand{\corollary}{\par \refstepcounter{corollary}%
{\bf Наслідок \arabic{corollary}.}}

\newcounter{remark}
\newcommand{\remark}{\par \refstepcounter{remark}%
{\bf Зауваження \arabic{remark}.}}

\newcounter{theorem}
\newcommand{\theorem}{\par \refstepcounter{theorem}%
{\bf Теорема \arabic{theorem}.}}

\newcounter{proposition}
\newcommand{\proposition}{\par \refstepcounter{proposition}%
{\bf Твердження \arabic{proposition}.}}

\newcounter{example}
\newcommand{\example}{\par \refstepcounter{example}%
{\bf Приклад \arabic{example}.}}

\renewcommand{\refname}{\centerline{\bf Список літератури}}

\renewcommand{\figurename}{Мал.}

\newcommand{\proof}{{\it Доведення.\,\,}}

\medskip\medskip
{\bf М.В.~Андрощук} (Житомирський державний університет імені Івана
Фран\-ка)

{\bf О.П.~Довгопятий} (Житомирський державний університет імені
Івана Фран\-ка)

{\bf Н.С.~Ількевич} (Житомирський державний
університет імені Івана Фран\-ка)

{\bf Є.О.~Севостьянов}
(Житомирський державний університет імені Івана Фран\-ка; Інститут
прикладної математики і механіки НАН України, м.~Слов'янськ)

\medskip\medskip\medskip\medskip
{\bf M.V.~Androschuk} (Zhytomyr Ivan Franko State University)

{\bf O.P.~Dovhopiatyi} (Zhytomyr Ivan Franko State University)

{\bf N.S.~Ilkevych} (Zhytomyr Ivan Franko State University)

{\bf E.A.~Sevost'yanov} (Zhytomyr Ivan Franko State University;
Institute of Applied Ma\-the\-ma\-tics and Mechanics of NAS of
Ukraine, Slov'yans'k)

\medskip\medskip\medskip
{\bf Про поведінку одного класу відображень, що діють на області з
локально квазіконформною межею}

{\bf On behavior of one class of mappings acting onto domains with a
locally quasiconformal boundary}

\medskip\medskip\medskip\medskip
Стаття присвячена дослідженню відображень, які задовольняють так
звану обернену нерівність Полецького. Розглядаються відображення
областей квазіекстремальної довжини, областей з локально
квазіконформною межею та регулярних (в сенсі простих кінців)
областей на області з локально квазіконформною межею, регулярних
областей, або областей, локально гельдерево еквівалентних до півкулі
на своїй межі. Для таких відображень нами отримана логарифмічна
неперервність за Гельдером в околі точок межі.

\medskip\medskip
The article is devoted to the study of mappings that satisfy the
so-called inverse Poletsky inequality. We consider mappings of
quasiextremal distance domains, domains with a locally
quasiconformal boundary, and domains which are regular in the sense
of prime ends onto domains with a locally quasiconformal boundary,
regular domains, or domains which are locally H\"{o}lder equivalent
to a half-ball on its boundary. For such mappings, we have obtained
the H\"{o}lder logarithmic continuity in some neighborhood of its
boundary points.

\newpage
{\bf 1. Вступ.} В одній нашій публікації ми розглянули відображення
одиничної кулі з оберненою умовою спотворення модуля сімей кривих
типу Полецького. Тут була встановлена їх логарифмічна неперервність
за Гельдером у межових точках (див.~\cite{Sev$_1$}). Дана замітка
присвячена розгляду цього питання в інших областях. Зокрема, ми
покажемо, що логарифмічна неперервність за Гельдером виконується в
межових точках заданої області, якщо ця область є областю
квазіекстремальної довжини, а відображена область є обмеженою
областю з локально квазіконформною межею. У статті розглянуті і інші
області, в тому числі, такі, відносно котрих логарифмічну
неперервність за Гельдером слід розуміти в сенсі простих кінців.
Логарифмічна неперервність за Гельдером у внутрішніх точках була
доведена раніше, причому у довільній області~\cite{SSD}. Відкрита
проблема стосується лише межових точок. Відзначимо також, що
обернені модульні нерівності відомі давно і відіграють ключову роль
при вивченні квазіконформних і квазірегулярних відображень, а також
відображень зі скінченним спотворенням довжини (див., напр.,
\cite[теорема~3.2]{MRV$_1$}, \cite[теорема~6.7.II]{Ri} та
\cite[теорема~8.5]{MRSY}).

\medskip
Нагадаємо деякі означення. Борелева функція $\rho:{\Bbb
R}^n\,\rightarrow [0,\infty] $ зветься {\it допустимою} щодо сім'ї
$\Gamma$ кривих $\gamma$ у ${\Bbb R}^n,$ якщо
\begin{equation}\label{eq1.4}
\int\limits_{\gamma}\rho (x)\, |dx|\geqslant 1
\end{equation}
для всіх (локально спрямованих) кривих $ \gamma \in \Gamma.$ У цьому
випадку ми пишемо: $\rho \in {\rm adm} \,\Gamma .$ {\it Модулем}
сім'ї кривих $\Gamma $ зветься величина
\begin{equation}\label{eq1.3gl0}
M(\Gamma)=\inf\limits_{\rho \in \,{\rm adm}\,\Gamma}
\int\limits_{{\Bbb R}^n} \rho^n (x)\,dm(x)\,.
\end{equation}
Нехай $Q:{\Bbb R}^n\rightarrow [0, \infty]$ -- вимірна за Лебегом
функція.  Будемо говорити, що {\it $f$ задовольняє обернену
нерівність Полецького}, якщо співвідношення
\begin{equation}\label{eq2*A}
M(\Gamma)\leqslant \int\limits_{f(D)} Q(y)\cdot\rho_*^n(y)\, dm(y)
\end{equation}
виконується для будь-якої сім'ї (локально спрямлюваних) кривих
$\Gamma$ в $D$ і для будь-якої $\rho_*\in {\rm adm}\,f(\Gamma).$
Зауважимо, що оцінки типу~(\ref{eq2*A}) добре відомі та виконуються
у багатьох класах відображень (див., напр.,
\cite[теорема~3.2]{MRV$_1$}, \cite[теорема~6.7.II]{Ri} та
\cite[теорема~8.5]{MRSY}).  Відображення $f:D\rightarrow {\Bbb R}^n$
називається {\it дискретним}, якщо прообраз
$\{f^{-1}\left(y\right)\}$ кожної точки $y\,\in\,{\Bbb R}^n$
складається з ізольованих точок, і {\it відкритим}, якщо образ
будь-якої відкритої множини $U\subset D$ є відкритою множиною в
${\Bbb R}^n.$ Відображення $f$ області $D$ на $D^{\,\prime}$
називається {\it замкненим}, якщо $f(E)$ є замкненим в
$D^{\,\prime}$ для будь-якої замкненої множини $E\subset D$ (див.,
напр., \cite[розд.~3]{Vu$_1$}).

\medskip
Область $D$ в ${\Bbb R}^n$ називається {\it областю
квазіекстремальної довжини} (скор. $QED$-{\it областю}), якщо
знайдеться таке число $A_0\geqslant 1,$ що для будь-яких континуумів
$E, F\subset D$ виконується нерівність
\begin{equation}\label{eq4***}
M(\Gamma(E, F, {\Bbb R}^n))\leqslant  A_0\cdot M(\Gamma(E, F, D))\,.
\end{equation}
Зауважимо, що одинична куля, півпростір або півкуля є областями
квазіекстремальної довжини, див.~\cite[лема~4.3]{Vu$_2$}.

\medskip
У подальшому,  в розширеному просторі $\overline{{{\Bbb
R}}^n}={{\Bbb R}}^n\cup\{\infty\}$ використовується {\it сферична
(хордальна) метрика} $h(x,y)=|\pi(x)-\pi(y)|,$ де $\pi$ --
стереографічна проекція $\overline{{{\Bbb R}}^n}$  на сферу
$S^n(\frac{1}{2}e_{n+1},\frac{1}{2})$ в ${{\Bbb R}}^{n+1},$ а саме,
$$h(x,\infty)=\frac{1}{\sqrt{1+{|x|}^2}}\,,$$
\begin{equation}\label{eq3C}
\ \ h(x,y)=\frac{|x-y|}{\sqrt{1+{|x|}^2} \sqrt{1+{|y|}^2}}\,, \ \
x\ne \infty\ne y
\end{equation}
(див., напр., \cite[означення~12.1]{Va}).
У подальшому, для множин $A, B\subset \overline{{\Bbb R}^n}$
покладемо
$$h(A, B)=\inf\limits_{x\in A, y\in B}h(x, y)\,,\quad h(A)=\sup\limits_{x, y\in A}h(x ,y)\,,$$
де $h$ -- хордальная відстань, визначена в~(\ref{eq3C}).

\medskip
Як звично, покладемо
$$
B(x_0, r)=\{x\in {\Bbb R}^n: |x-x_0|<r\}\,,\qquad {\Bbb B}^n=B(0,
1)\,,$$
$$S(x_0,r) = \{ x\,\in\,{\Bbb R}^n : |x-x_0|=r\}\,.$$
Розглянемо наступне означення, яке було запропоновано
Няккі~\cite{Na}, див. також~\cite{KR$_1$}. Межа області $D$ в ${\Bbb
R}^n$ називається {\it локально квазіконформною,} якщо кожна точка
$x_0\in\partial D$ має окіл $U,$ для якого існує квазіконформне
відображення $\varphi$ околу $U$ на одиничну кулю ${\Bbb
B}^n\subset{\Bbb R}^n$ таке, що $\varphi(\partial D\cap U)$ є
перетином одиничної кулі ${\Bbb B}^n$ з координатною гіперплощиною
$x_n=0,$ де $x=(x_1,\ldots, x_n.)$

\medskip
Для числа $\delta>0,$ областей $D, D^{\,\prime}\subset {\Bbb R}^n,$
$n\geqslant 2,$ невиродженого континуума $A\subset D^{\,\prime}$ і
вимірної за Лебегом функції $Q:D^{\,\prime}\rightarrow [0, \infty]$
позначимо через ${\frak S}_{\delta, A, Q }(D, D^{\,\prime})$ сім'ю
всіх відкритих дискретних і замкнених відображень $f$ області $D$ на
область $D^{\,\prime},$ що задовольняють умову~(\ref{eq2*A}) і
таких, що $h(f^{\,-1}(A),
\partial D)\geqslant~\delta.$ Справедливо наступне твердження.

\medskip
\begin{theorem}\label{th1}
{\sl\,Нехай $Q\in L^1(D^{\,\prime}),$ $D$ є областю
квазіекстремальної довжини, а $D^{\,\prime}$ є обмеженою областю з
локально квазіконформною межею. Тоді будь-яке відображення $f\in
{\frak S}_{\delta, A, Q}(D, D^{\,\prime}),$ яке задовольняє
співвідношення~(\ref{eq2*A}), продовжується до відображення
$f:\overline{D}\rightarrow\overline{D^{\,\prime}},$ при цьому, для
кожної точки $x_0\in\partial D$ знайдуться окіл $U$ цієї точки і
сталі $C_n=C(n, A, D, D^{\,\prime}, x_0)>0$ і $0<\alpha=\alpha(n, A,
D, D^{\,\prime}, x_0)\leqslant 1$ така, що
\begin{equation}\label{eq2C}
|\overline{f}(x)-\overline{f}(y)|^{\frac{n}{\alpha^2}}\leqslant\frac{C_n\cdot
\Vert Q\Vert_1}{\log\left(1+\frac{\delta}{2|x-y|}\right)}
\end{equation}
для всіх $x, y\in U\cap \overline{D},$
де $\Vert Q\Vert_1$ -- норма функції $Q$ в $L^1(D).$

 }
\end{theorem}

\medskip
{\bf 2. Допоміжні леми.} Наступна лема у випадку одиничної кулі була
доведена в~\cite[лема~2.1]{Sev$_1$}, а в випадку довільної обмеженої
опуклої області -- в~\cite[лема~2.1]{DS}.

\medskip
\begin{lemma}\label{lem1}
{\sl\, Нехай $E$ -- довільний континуум, що належить
області~$D^{\,\prime},$ $Q\in L^1(D^{\,\prime}).$ Тоді існує
$\delta_1>0$ таке, що ${\frak S}_{\delta, A, Q }\subset {\frak
S}_{\delta_1, E, Q}.$ Іншими словами, якщо $f$ -- відкрите дискретне
і замкнене відображення області $D$ на $D^{\,\prime}$ з
умовою~(\ref{eq2*A}), таке, що $h(f^{\,-1}(A),
\partial D)\geqslant~\delta,$ то існує $\delta_1>0,$ не залежне від
$f$ таке, що $h(f^{\,-1}(E),
\partial D)\geqslant~\delta_1.$ }
\end{lemma}

\medskip
Наступна лема була доведена в випадку, коли область $D^{\,\prime}$ є
одиничною кулею (див. хід доведення теореми~1.1 в~\cite{Sev$_1$}).
Для довільної обмеженої опуклої області доведення здійснено
в~\cite[лема~2.2]{DS}.

\medskip
\begin{lemma}\label{lem2}
{\sl\, Нехай $D_1$ -- обмежена опукла область в ${\Bbb R}^n,$
$n\geqslant 2,$ і нехай $B(y_*, \delta_*/2)$ -- куля з центром в
точці $y_*\in D_1,$ де $\delta_*:=d(y_*, \partial D_1).$ Нехай
$z_0\in
\partial D_1.$ Тоді для будь-яких точок $A, B\in B(z_0,
\delta_*/8)\cap D_1$ знайдуться точки $C, D\in \overline{B(y_*,
\delta_*/2)},$ для яких відрізки $[A, C]$ і $[B, D]$ є такими, що
\begin{equation}\label{eq13}
{\rm dist\,}([A, C], [B, D])\geqslant C_0\cdot |A-B|\,,
\end{equation}
де $C_0>0$ -- деяка стала, яка залежить тільки від $\delta_*$ і
$d(D_1).$
 }
\end{lemma}

\medskip
{\bf 3. Відображення областей квазіекстремальної довжини. Доведення
теореми~\ref{th1}}. Можливість неперервного продовження відображення
$f$ на межу області $D$ випливає з теореми~3.1 в~\cite{SSD}.
Зокрема, слабка плоскість $\partial D$ є наслідком того, що $D$ є
$QED$-областю (див., напр., \cite[лема~2]{SevSkv}), а те, що область
з локально квазіконформною межею є локально зв'язною на своїй межі є
наслідком означення цієї області.

\medskip
Доведемо логарифмічну неперервність за Гельдером~(\ref{eq2C}).
Достатньо довести співвідношення~(\ref{eq2C}) для $x, y\in U\cap D.$

Будемо міркувати методом від супротивного, а саме, припустимо, що
співвідношення~(\ref{eq2C}) не виконується принаймні в одній
точці~$x_0\in \partial D.$ Тоді для будь-якого натурального $m\in
{\Bbb N}$ знайдуться точки $x_m, y_m\in D$ і відображення $f_m\in
{\frak S}_{\delta, A, Q}(D, D^{\,\prime}),$ такі що
$|x_m-x_0|<\frac{1}{m},$ $|y_m-x_0|<\frac{1}{m},$ проте
\begin{equation}\label{eq2D}
{|f_m(x_m)-f_m(y_m)|}^{m}> m\cdot  \frac{\Vert
Q\Vert_1}{\log\left(1+\frac{\delta}{2|x_m-y_m|}\right)}\,.
\end{equation}
Оскільки область $D^{\,\prime}$ обмежена, ми можемо вважати, що
\begin{equation}\label{eq1}
|f_m(x_m)-f_m(y_m)|<1\,, \quad |x_m-y_m|<1\,,\quad  m=1,2,\ldots \,,
\end{equation}
і що послідовності $f_m(x_m)$ і $f_m(y_m)$ збігаються при
$m\rightarrow\infty.$

Зауважимо, що в цьому випадку $f_m(x_m)$ і $f_m(y_m)$ є збіжними при
$m\rightarrow\infty$ до однієї і тієї ж самої точки $y_0\in \partial
D^{\,\prime}.$ Дійсно, в силу компактності $\partial D^{\,\prime}$
послідовність $f_m(x_0),$ $m=1,2,\ldots,$ також можна вважати
збіжною при $m\rightarrow\infty,$ причому, ця послідовність
збігається до деякої точки з $\partial D^{\,\prime},$ бо межа будь
якої області є замкненою. Нехай $f_m(x_0)\rightarrow y_0$ при
$m\rightarrow\infty.$

Зауважимо, що за теоремою~1.2 в~\cite{SSD} сім'я відображень~${\frak
S}_{\delta, A, Q}$ є одностайно неперервною в~$\overline{D}.$ Тоді
за нерівністю трикутника
$$|f_m(x_m)-y_0|\leqslant |f_m(x_m)-f_m(x_0)|+|f_m(x_0)-y_0|\rightarrow 0\,,\quad m\rightarrow\infty
\,.$$
Останнє співвідношення доводить, що $f_m(x_m)\rightarrow y_0\in
\partial D^{\,\prime}$ при $m\rightarrow\infty.$ Аналогічно можна
довести, що $f_m(y_m)\rightarrow y_0\in
\partial D^{\,\prime}$ при $m\rightarrow\infty.$

\medskip
Скористаємося тим, що межа області $D^{\,\prime}$ є локально
квазіконформною. Нехай $U_0$ -- окіл точки $y_0,$ для якого існує
квазіконформне відображення $\varphi:U_0\rightarrow {\Bbb B}^n$
таке, що $\varphi(U_0)={\Bbb B}^n,$ $\varphi(U_0\cap
D^{\,\prime})={\Bbb B}^n_+,$
\begin{equation}\label{eq24}
{\Bbb B}^n_+=\{x=(x_1,\ldots x_n)\in {\Bbb B}^n: x_n>0\}\,.
\end{equation}
Застосовуючи додаткове мебіусове перетворення можна вважати, що
$\varphi(y_0)=0$ (див.~\cite[доведення теореми~17.10]{Va}).

Покладемо $z_0=0,$ $y_*:=\left(0,\ldots, 0, \frac{1}{2}\right),$
$D_1:={\Bbb B}^n_+.$ Тоді $\delta_*:={\rm dist}\,(y_*, \partial
{\Bbb B}^n_+)=\frac{1}{2}.$ Оскільки $f_m(x_m)\rightarrow y_0\in
\partial D^{\,\prime}$ і $f_m(y_m)\rightarrow y_0\in
\partial D^{\,\prime}$ при $m\rightarrow\infty,$ то $\varphi(f_m(x_m))\rightarrow 0\in
\partial D^{\,\prime}$ і $\varphi(f_m(y_m))\rightarrow 0\in
\partial D^{\,\prime}$ при $m\rightarrow\infty.$ Тоді
$\varphi(f_m(x_m))$ і $\varphi(f_m(y_m))\in B(0, 1/16)$ при всіх
$m>m_1$ і деякому $m_1\in {\Bbb N}.$

З огляду на опуклість ${\Bbb B}^n_+$ можна застосувати
лему~\ref{lem2}. За цією лемою знайдуться точки $C_m, D_m\in
\overline{B(y_*, 1/4)},$ для яких відрізки $[\varphi(f_m(x_m)),
C_m]$ і $[\varphi(f_m(y_m)), D_m]$ є такими, що
\begin{equation}\label{eq13A}
{\rm dist\,}([\varphi(f_m(x_m)), C_m], [\varphi(f_m(y_m)),
D_m])\geqslant C_0\cdot
|\varphi(f_m(x_m))-\varphi(f_m(y_m))|\,,\quad m\geqslant m_1\,,
\end{equation}
де $C_0>0$ -- деяка абсолютна стала. Покладемо
$$E_*:=\varphi^{\,-1}(\overline{B(y_*, 1/4)})\,,$$
%
$$\alpha_m:=\varphi^{\,-1}([\varphi(f_m(x_m)), C_m]),
\quad \beta_m:=\varphi^{\,-1}([\varphi(f_m(x_m)), D_m])\,.$$
Нехай
$${\rm dist\,}(|\alpha_m|, |\beta_m|)=|w_m-u_m|\,,\qquad w_m\in |\alpha_m|, u_m\in |\beta_m|\,.$$
З огляду на те, що квазіконформні відображення $\varphi$ та
$\varphi^{\,-1}$ є локально гельдеровими з деяким показником
$0<\alpha\leqslant 1$ та сталими гельдеровості $C_1>0$ та
$C^{\,*}_1,$ відповідно (див.~\cite[теорема~1.11.III]{Ri}), а також
з огляду на~(\ref{eq13A}), ми будемо мати, що
$${\rm dist\,}(|\alpha_m|, |\beta_m|)=|w_m-u_m|\geqslant
\frac{1}{C^{\frac{1}{\alpha}}_1}\cdot
|\varphi(w_m)-\varphi(u_m)|^{\frac{1}{\alpha}}\geqslant $$
\begin{equation}\label{eq22}
\geqslant \frac{1}{C^{\frac{1}{\alpha}}_1}\cdot \left({\rm
dist\,}([\varphi(f_m(x_m)), C_m], [\varphi(f_m(y_m)),
D_m])\right)^{\frac{1}{\alpha}}\geqslant
\frac{C^{\frac{1}{\alpha}}_0}{C^{\frac{1}{\alpha}}_1}\cdot
(|\varphi(f_m(x_m))-\varphi(f_m(y_m))|)^{\frac{1}{\alpha}}\geqslant
\end{equation}
$$\geqslant \frac{C^{\frac{1}{\alpha}}_0}
{C^{\frac{1}{\alpha}}_1\left(C^{\,*}_1\right)^{\frac{1}{\alpha^2}}}\cdot
|f_m(x_m)-f_m(y_m)|^{\frac{1}{\alpha^2}}\,.$$

\medskip
Нехай $\alpha^{\,*}_m,$ $\beta^{\,*}_m$ -- повні $f_m$-підняття
кривих $\alpha_m$ і $\beta_m$ з початками в точках $x_m$ і $y_m,$
відповідно (вони існують за~\cite[лема~3.7]{Vu$_1$}).
Тоді за означенням $|\alpha^{\,*}_m|\cap
f_m^{\,-1}(E)\ne\varnothing\ne |\beta^{\,*}_m|\cap f^{\,-1}(E).$ За
лемою~\ref{lem1} існує $\delta_1>0$ таке, що $h(f_m^{\,-1}(E),
\partial D)\geqslant \delta_1$ для всіх $m\in {\Bbb N}.$
Оскільки $x_m, y_m\in B(x_0, 1/m),$ то для достатньо великих $m\in
{\Bbb N}$
\begin{equation}\label{eq4}
d(|\alpha^{\,*}_m|)\geqslant \delta_1/2\,,\quad
d(|\beta^{\,*}_m|)\geqslant \delta_1/2\,.
\end{equation}
Нехай
$$\Gamma_m:=\Gamma(|\alpha^{\,*}_m|, |\beta^{\,*}_m|, D)\,.$$
Тоді з одного боку за нерівністю~(\ref{eq4***})
\begin{equation}\label{eq7A}
M(\Gamma_m)\geqslant (1/A_0)\cdot M(\Gamma_m(|\alpha^{\,*}_m|,
|\beta^{\,*}_m|, {\Bbb R}^n))\,,
\end{equation}
а з іншого боку, за~\cite[лема~7.38]{Vu$_3$}
\begin{equation}\label{eq7B}
M(\Gamma_m(|\alpha^{\,*}_m|, |\beta^{\,*}_m|, {\Bbb R}^n))\geqslant
c_n\cdot\log\left(1+\frac{1}{\widetilde{m}}\right)\,,
\end{equation}
де $c_n>0$ -- деяка стала, яка залежить лише від $n,$
$$\widetilde{m}=\frac{{\rm dist}(|\alpha^{\,*}_m|,
|\beta^{\,*}_m|)}{\min\{{\rm diam\,}(|\alpha^{\,*}_m)|, {\rm
diam\,}(|\beta^{\,*}_m|)\}}\,.$$
Тоді поєднуючи~(\ref{eq4}), (\ref{eq7A}) і~(\ref{eq7B}) і
враховуючи, що ${\rm dist}\,(|\alpha^{\,*}_m|,
|\beta^{\,*}_m|)\leqslant |x_m-y_m|,$  ми отримуємо, що
\begin{equation}\label{eq7C}
M(\Gamma_m)\geqslant \widetilde{c_n}\cdot
\log\left(1+\frac{\delta_1}{2{\rm dist}(|\alpha^{\,*}_m|,
|\beta^{\,*}_m|)}\right)\geqslant \widetilde{c_n}\cdot
\log\left(1+\frac{\delta_1}{2|x_m-y_m|}\right)\,,
\end{equation}
де $\widetilde{c_n}>0$ -- деяка стала, яка залежить тільки від $n$ і
сталої $A_0$ з означення $QED$-області.

\medskip
Встановимо тепер верхню оцінку для $M(\Gamma_m).$ Покладемо
$$\rho_m(y)= \left\{
\begin{array}{rr}
\frac{C^{\frac{1}{\alpha}}_1\left(C^{\,*}_1\right)^{\frac{1}{\alpha^2}}}{C^{\frac{1}{\alpha}}_0}\cdot
|f_m(x_m)-f_m(y_m)|^{-\frac{1}{\alpha^2}}, & y\in D^{\,\prime},\\
0,  &  y\not\in D^{\,\prime}\,.
\end{array}
\right. $$
Зауважимо, що $\rho_m$ задовольняє співвідношення~(\ref{eq1.4}) для
сім'ї кривих $f_m(\Gamma_m)$ в силу співвідношення~(\ref{eq22}).
Тоді за означення сім'ї ${\frak S}_{\delta, A, Q }$ ми отримаємо, що
\begin{equation}\label{eq14***}
M(\Gamma_m)\leqslant
\frac{\left(\frac{C^{\frac{1}{\alpha}}_1\left(C^{\,*}_1\right)^{\frac{1}{\alpha^2}}}{C^{\frac{1}{\alpha}}_0}\right)^n}
{|f_m(x_m)-f_m(y_m)|^{\frac{n}{\alpha^2}}} \cdot
\int\limits_{D^{\,\prime}} Q(y)\,dm(y)\,.
\end{equation}
З~(\ref{eq7C}) і (\ref{eq14***}) випливає, що
$$\widetilde{c_n}\cdot \log\left(1+\frac{\delta_1}{2|x_m-y_m|}\right)\leqslant
\left(\frac{C^{\frac{1}{\alpha}}_1\left(C^{\,*}_1\right)^{\frac{1}{\alpha^2}}}{C^{\frac{1}{\alpha}}_0}\right)^n
\cdot\frac{\Vert
Q\Vert_1}{{|f_m(x_m)-f_m(y_m)|}^{\frac{n}{\alpha^2}}}\,.$$
З останнього співвідношення випливає, що
\begin{equation}\label{eq3D}
|f_m(x_m)-f_m(y_m)|^{\frac{n}{\alpha^2}}\leqslant
(\widetilde{c_n})^{-1}\cdot
\left(\frac{C^{\frac{1}{\alpha}}_1\left(C^{\,*}_1\right)^{\frac{1}{\alpha^2}}}{C^{\frac{1}{\alpha}}_0}\right)^{n}
\cdot\frac{\left(\Vert Q\Vert_1\right)}
{\log\left(1+\frac{\delta_1}{2|x_m-y_m|}\right)}\,.
\end{equation}
Покладемо $$M_0:=\max\left\{\left[(\widetilde{c_n})^{-1}\cdot
\left(\frac{C^{\frac{1}{\alpha}}_1\left(C^{\,*}_1\right)^{\frac{1}{\alpha^2}}}{C^{\frac{1}{\alpha}}_0}\right)^{n}
\right]+1, \left[\frac{n}{\alpha^2}\right]+1\right\}\,.$$
З огляду на співвідношення~(\ref{eq1}),
$|f_m(x_m)-f_m(y_m)|^{M_0}\leqslant
|f_m(x_m)-f_m(y_m)|^{\frac{n}{\alpha^2}}$ для достатньо великих
$m>m_1.$ Тоді з~(\ref{eq3D}) випливає, що
\begin{equation}\label{eq3DA}
|f_m(x_m)-f_m(y_m)|^{M_0}\leqslant M_0 \cdot\frac{\Vert Q\Vert_1}
{\log\left(1+\frac{\delta_1}{2|x_m-y_m|}\right)}\,,\quad m>m_1.
\end{equation}
Нарешті, з огляду на правило Лопіталя,
$\log\left(1+\frac{1}{nt}\right)\sim\log\left(1+\frac{1}{kt}\right)$
при $t\rightarrow+0$ для різних фіксованих $k, n> 0.$ Тоді при
деякій натуральній сталій $M_1\geqslant M_0$ виконується
співвідношення $\frac{\Vert Q\Vert_1}
{\log\left(1+\frac{\delta_1}{2|x_m-y_m|}\right)}\leqslant
\frac{M_1}{M_0}\cdot\frac{\Vert Q\Vert_1}
{\log\left(1+\frac{\delta}{2|x_m-y_m|}\right)},$ $m>m_2>m_1.$ Тоді
з~(\ref{eq3DA}) ми отримаємо, що
\begin{equation}\label{eq4D}
|f_m(x_m)-f_m(y_m)|^{M_1}\leqslant M_1 \cdot\frac{\Vert Q\Vert_1}
{\log\left(1+\frac{\delta}{2|x_m-y_m|}\right)}\,,\quad m>m_2\,.
\end{equation}
Тоді з~(\ref{eq4D}) при $m>N_0:=\max\{m_2, [M_1]+1\}$ ми будемо
мати, що
\begin{equation}\label{eq5D}
|f_m(x_m)-f_m(y_m)|^{m}\leqslant m \cdot\frac{\Vert Q\Vert_1}
{\log\left(1+\frac{\delta}{2|x_m-y_m|}\right)}\,,\quad m>N_0\,.
\end{equation}
Співвідношення~(\ref{eq5D}) суперечить припущенню, зробленому
в~(\ref{eq2D}). Отримана суперечність доводить теорему.~$\Box$

\medskip
{\bf 4. Відображення на області, локально-гельдерево та
локально-ліпшицево еквівалентні до півкулі на межі.} По аналогії з
областями, які мають локально квазіконформні межі, впровадимо
наступні означення, які допоможуть нам сформулювати результати щодо
логарифмічної неперервності за Гельдером з конкретними показниками
$\alpha.$

\medskip
Домовимось називати область $D$ в ${\Bbb R}^n$ {\it локально
$\alpha$-гельдерово еквівалентною до півкулі на своїй межі,} якщо
для кожної точки $x_0\in\partial D$ існує її окіл $U,$ сталі
$C_1=C_1(x_0)>0,$ $C^{\,*}_1=C^{\,*}_1(x_0)>0$ і гомеоморфізм
$\varphi$ околу $U$ на ${\Bbb B}^n$ такий що
\begin{equation}\label{eq25}
\frac{1}{\left(C^{\,*}_1\right)^{\frac{1}{\alpha}}}|x-y|^{\frac{1}{\alpha}}\leqslant
|\varphi(x)-\varphi(y)|\leqslant C\cdot
|x-y|^{\alpha}\qquad\forall\,\,x,y\in U\,,
\end{equation}
причому $\varphi(\partial D\cap U)={\Bbb B}^n_+,$ де ${\Bbb B}^n_+$
визначено в~(\ref{eq24}). Область $D$ в ${\Bbb R}^n$ буде називатися
{\it локально ліпшицево еквівалентною до півкулі на своїй межі,}
якщо в~(\ref{eq25}) $\alpha=1.$ Є справедливим наступне твердження.

\medskip
\begin{theorem}\label{th5}
{\sl\,Нехай $Q\in L^1(D^{\,\prime}),$ $D$ є областю
квазіекстремальної довжини, а $D^{\,\prime}$ є обмеженою областю,
яка є локально $\alpha$-гельдерово еквівалентною до півкулі на своїй
межі. Тоді будь-яке відображення $f\in {\frak S}_{\delta, A, Q}(D,
D^{\,\prime}),$ яке задовольняє співвідношення~(\ref{eq2*A}),
продовжується до відображення
$f:\overline{D}\rightarrow\overline{D^{\,\prime}},$ при цьому, для
кожної точки $x_0\in\partial D$ знайдуться окіл $U$ цієї точки і
стала $C_n=C(n, A, D, D^{\,\prime}, x_0)>0$ така, що
\begin{equation}\label{eq23}
|\overline{f}(x)-\overline{f}(y)|^{\frac{n}{\alpha^2}}\leqslant\frac{C_n\cdot
\Vert Q\Vert_1}{\log\left(1+\frac{\delta}{2|x-y|}\right)}
\end{equation}
для всіх $x, y\in U\cap \overline{D},$
де $\Vert Q\Vert_1$ -- норма функції $Q$ в $L^1(D).$

 }
\end{theorem}

\medskip
\begin{proof}
Хід доведення теореми~\ref{th5} майже дослівно повторює хід
доведення попередньої теореми~\ref{th1}, тому обмежимося лише
схематичним доведенням. Достатньо довести
співвідношення~(\ref{eq23}) для $x, y\in U\cap D$ (випадок $x, y\in
U\cap \overline{D}$ тоді може бути отриманий шляхом операції
граничного переходу).

\medskip
Як і при доведенні попередньої теореми~\ref{th1}, будемо міркувати
методом від супротивного. Припустимо, що співвідношення~(\ref{eq23})
не виконується принаймні в одній точці~$x_0\in \partial D.$ Тоді для
будь-якого натурального $m\in {\Bbb N}$ знайдуться точки $x_m,
y_m\in D$ і відображення $f_m\in {\frak S}_{\delta, A, Q}(D,
D^{\,\prime}),$ такі що $|x_m-x_0|<\frac{1}{m},$
$|y_m-x_0|<\frac{1}{m},$ проте
\begin{equation}\label{eq26}
{|f_m(x_m)-f_m(y_m)|}^{\frac{n}{\alpha^2}}> m\cdot  \frac{\Vert
Q\Vert_1}{\log\left(1+\frac{\delta}{2|x_m-y_m|}\right)}\,.
\end{equation}
Оскільки область $D^{\,\prime}$ обмежена, ми можемо вважати, що
\begin{equation}\label{eq27}
|f_m(x_m)-f_m(y_m)|<1\,, \quad |x_m-y_m|<1\,,\quad  m=1,2,\ldots \,,
\end{equation}
і що послідовності $f_m(x_m)$ і $f_m(y_m)$ збігаються при
$m\rightarrow\infty.$

\medskip
Далі, використовуючи позначення з доведення теореми~\ref{th1},
приходимо до співвідношення~(\ref{eq22}), яке виходить з означення
відображення~$\varphi$ у~(\ref{eq25}). Міркуючи далі аналогічно
доведенню теореми~\ref{th1}, приходимо до
співвідношення~(\ref{eq3DA}), яке суперечить припущенню~(\ref{eq26})
(тут слід врахувати, що за правилом Лопіталя
$\log\left(1+\frac{1}{nt}\right)\sim\log\left(1+\frac{1}{kt}\right)$
при $t\rightarrow+0$ для різних фіксованих $k, n> 0$). Отримана
суперечність завершує доведення теореми~(\ref{th5}).~$\Box$
\end{proof}

\medskip
{\bf 5. Відображення областей з квазіконформними межами на такі ж
області.}  Справедлива наступна

\medskip
\begin{theorem}\label{th3}
{\sl\, Нехай $Q\in L^1(D^{\,\prime}),$ $D$ і $D^{\,\prime}$ --
області з локально квазіконформною межею. Тоді будь-яке відображення
$f\in {\frak S}_{\delta, A, Q}(D, D^{\,\prime}),$ яке задовольняє
співвідношення~(\ref{eq2*A}), продовжується до відображення
$f:\overline{D}\rightarrow\overline{D^{\,\prime}},$ при цьому, для
кожної точки $x_0\in\partial D$ знайдуться окіл $U$ цієї точки і
сталі $C_n=C(n, A, D, D^{\,\prime}, x_0)>0,$
$0<\alpha^{\,*}=\alpha^{\,*}(n, A, D, D^{\,\prime}, x_0)\leqslant 1$
і $0<\alpha=\alpha(n, A, D, D^{\,\prime}, x_0)\leqslant 1$ така, що
\begin{equation}\label{eq28}
|\overline{f}(x)-\overline{f}(y)|^{\frac{n}{\alpha^2}}\leqslant\frac{C_n\cdot
\Vert
Q\Vert_1}{\log\left(1+\frac{\delta}{2|x-y|^{\,\alpha^{\,*}}}\right)}
\end{equation}
для всіх $x, y\in U\cap \overline{D},$
де $\Vert Q\Vert_1$ -- норма функції $Q$ в $L^1(D).$

 }
\end{theorem}

\medskip
\begin{proof} Можливість неперервного продовження відображення $f$ на межу області
$D$ випливає з теореми~3.1 в~\cite{SSD}. Зокрема, області з локально
квазіконформними межами є слабко плоскими (див.
\cite[Proposition~2.2]{ISS}, або~\cite[теорема~17.10]{Va}). З іншого
боку, очевидно, такі області є й локально зв'язними на своїй межі.

\medskip
За означенням локально квазіконформної межі області $D$ існує окіл
$U^{\,*}$ точки $x_0$ і квазіконформне відображення
$\varphi^{\,*}:U^{\,*}\rightarrow {\Bbb B}^n,$
$\varphi(U^{\,*})={\Bbb B}^n,$ таке що $\varphi^{\,*}(D\cap
U^{\,*})={\Bbb B}_{+}^n,$ де ${\Bbb B}_{+}^n=\{x\in {\Bbb B}^n:
x=(x_1,\ldots, x_n), x_n>0\}$ -- півкуля. Можна вважати, що $x_0\ne
\infty$ і $\varphi^{\,*}(x_0)=0$ (див. хід доведення теореми~17.10
у~\cite{Va}). Оскільки $\varphi^{\,*}$ є квазіконформним
відображенням, таким є і відображення $(\varphi^{\,*})^{\,-1},$
причому воно є неперервним за Гельдером з деякою сталою
$\widetilde{C}>0$ і деяким показником $0<\alpha^{\,*}\leqslant 1$
(див.~\cite[теорема~1.11.III]{Ri}).

\medskip
Доведемо співвідношення~(\ref{eq28}) зі вказаним
показником~$\alpha^{\,*}.$ Як і при доведенні попередніх теорем,
використаємо метод від супротивного. Припустимо, що
співвідношення~(\ref{eq28}) не виконується принаймні в одній
точці~$x_0\in \partial D.$ Тоді для будь-якого натурального $m\in
{\Bbb N}$ знайдуться точки $x_m, y_m\in D$ і відображення $f_m\in
{\frak S}_{\delta, A, Q}(D, D^{\,\prime}),$ такі що
$|x_m-x_0|<\frac{1}{m},$ $|y_m-x_0|<\frac{1}{m},$ проте виконується
співвідношення
\begin{equation}\label{eq2E}
{|f_m(x_m)-f_m(y_m)|}^{m}> m\cdot  \frac{\Vert
Q\Vert_1}{\log\left(1+\frac{\delta}{2|x_m-y_m|^{\alpha^{\,*}}}\right)}\,.
\end{equation}

\medskip
Далі, міркуючи аналогічно доведенню теореми~\ref{th1}, ми отримаємо
співвідношення~(\ref{eq22}) (в якому ми зберігаємо всі позначення,
впровадженні під час цього доведення).

\medskip Нехай $V^{\,*}$ -- окіл точки $x_0,$ який належить $U^{\,*}$
разом із своїм замиканням. Покладемо
\begin{equation}\label{eq23A}
\delta_2:={\rm dist\,}(\partial V^{\,*}, \partial U^{\,*})\,.
\end{equation}
Без обмеження загальності, можна вважати, що $x_m, y_m\in V^{\,*}$
при всіх $m\in {\Bbb N}.$ Розглянемо допоміжне відображення
\begin{equation}\label{eq16}
F_m(w):=f_m((\varphi^{\,*})^{\,-1}(w))\,,\qquad F_m:{\Bbb
B}_{+}^n\rightarrow D^{\,\prime}\,.
\end{equation}
Нехай $\alpha^{\,*}_m,$ $\beta^{\,*}_m$ -- повні $f_m$-підняття
кривих $\alpha_m$ і $\beta_m$ з початками в точках $x_m$ і $y_m,$
відповідно (вони існують за~\cite[лема~3.7]{Vu$_1$}). Тоді за
означенням $|\alpha^{\,*}_m|\cap f^{\,-1}(E)\ne\varnothing\ne
|\beta^{\,*}_m|\cap f^{\,-1}(E).$ Оскільки $d(f_m^{\,-1}(E),
\partial D)\geqslant \delta_1$ і $x_m, y_m\in V,$ то
$$|\alpha^{\,*}_m|\cap U^{\*}\ne\varnothing\ne |\alpha^{\,*}_m|\cap({\Bbb R}^n\setminus
U^{\*})$$ і
$$|\beta^{\,*}_m|\cap U^{\*}\ne\varnothing\ne |\beta^{\,*}_m|\cap({\Bbb R}^n\setminus
U^{\*})\,.$$
З огляду на~\cite[теорема~1.I.5.46]{Ku}
\begin{equation}\label{eq31} |\alpha^{\,*}_m|\cap \partial U^{\*}\ne\varnothing\,,
|\beta^{\,*}_m|\cap
\partial U^{\*}\ne\varnothing\,.
\end{equation}
Аналогічно,
\begin{equation}\label{eq32} |\alpha^{\,*}_m|\cap \partial V\ne\varnothing\,,
|\beta^{\,*}_m|\cap
\partial V\ne\varnothing\,.
\end{equation}
З огляду на~(\ref{eq32}), криві $\alpha^{\,*}_m$ і $\beta^{\,*}_m$
мають підкриві $\alpha^{\,*}_m$ і $\beta^{\,*}_m$ з початками в
точках $x_m$ і $y_m,$ які мають початки у $V^{\,*}$ та мають кінцеві
точки у $\partial U^{\,*}.$ За співвідношеннями~(\ref{eq23A}),
(\ref{eq31}) і (\ref{eq32}) ми будемо мати, що
\begin{equation}\label{eq4_1}
d(\alpha^{\,*}_m)\geqslant \delta_2\,,\quad
d(\beta^{\,*}_m)\geqslant \delta_2\,.
\end{equation}
Розглянемо криві $\varphi^{\,*}(\alpha^{\,*}_m)$ і
$\varphi^{\,*}(\beta^{\,*}_m).$
Нехай $\overline{x}_m, \overline{y}_m\in U^{\,*}$ є такими, що
$d(\alpha^{\,*}_m)=|\overline{x}_m-\overline{y}_m|.$ Покладемо
$x^{\,*}_m=\varphi^{\,*}(\overline{x}_m)$ і
$y^{\,*}_m=\varphi^{\,*}(\overline{y}_m).$ Тоді
$$|x^{\,*}_m-y^{\,*}_m|^{\alpha^{\,*}}\geqslant\frac{1}{\widetilde{C}}\cdot
|\overline{x}_m-\overline{y}_m|=d(\alpha^{\,*}_m)\geqslant
\frac{1}{\widetilde{C}}\delta_2\,,$$
або
\begin{equation}\label{eq15}
|x^{\,*}_m-y^{\,*}_m|\geqslant
\left(\frac{1}{\widetilde{C}}\delta_2\right)^{1/\alpha^{\,*}}\,.
\end{equation}
З~(\ref{eq15}) отримаємо, що
$d(\varphi^{\,*}(\alpha^{\,*}_m))\geqslant
\left(\frac{1}{\widetilde{C}}\delta_2\right)^{1/\alpha^{\,*}}.$
Аналогічно, $d(\varphi^{\,*}(\beta^{\,*}_m))\geqslant
\left(\frac{1}{\widetilde{C}}\delta_2\right)^{1/^{\alpha^{\,*}}}.$
Нехай
$$\Gamma_m:=\Gamma(\varphi^{\,*}(\alpha^{\,*}_m), \varphi^{\,*}(\beta^{\,*}_m), {\Bbb B}_+^n)\,.$$
Зауважимо, що~${\Bbb B}_+^n$ є обмеженою опуклою областю, тому вона
є областю Джона (див.~\cite[зауваження~2.4]{MS}), отже, є
рівномірною областю (див.~\cite[зауваження~2.13(c)]{MS}), тому є
також і $QED$-областю з деяким $A_0^{\,*}<\infty$ у~(\ref{eq4***})
(див.~\cite[лема~2.18]{GM}). Тоді з одного боку за
нерівністю~(\ref{eq4***})
\begin{equation}\label{eq7A_1}
M(\Gamma_m)\geqslant (1/A^{\,*}_0)\cdot
M(\Gamma(\varphi^{\,*}(\alpha^{\,*}_m),
\varphi^{\,*}(\beta^{\,*}_m), {\Bbb R}^n))\,,
\end{equation}
а з іншого боку, за~\cite[лема~7.38]{Vu$_3$}
\begin{equation}\label{eq7B_1}
M(\Gamma(\varphi^{\,*}(\alpha^{\,*}_m),
\varphi^{\,*}(\beta^{\,*}_m), {\Bbb R}^n))\geqslant
c_n\cdot\log\left(1+\frac{1}{\widetilde{m}}\right)\,,
\end{equation}
де $c_n>0$ -- деяка стала, яка залежить лише від $n,$
$$\widetilde{m}=\frac{{\rm dist}(\varphi^{\,*}(\alpha^{\,*}_m), \varphi^{\,*}(\alpha^{\,*}_m))}
{\min\{{\rm diam\,}(\varphi^{\,*}(\alpha^{\,*}_m)), {\rm
diam\,}\varphi^{\,*}(\beta^{\,*}_m)\}}\,.$$
Тоді поєднуючи~(\ref{eq7A_1}) і~(\ref{eq7B_1}) і враховуючи, що
${\rm dist}\,(\varphi^{\,*}(\alpha^{\,*}_m),
\varphi^{\,*}(\beta^{\,*}_m))\leqslant
|\varphi^{\,*}(x_m)-\varphi^{\,*}(y_m)|,$ ми отримуємо, що
$$M(\Gamma_m)\geqslant \widetilde{c_n}\cdot \log\left
(1+\frac{\delta_2^{1/\alpha^{\,*}}}{(\widetilde{C})^{1/\alpha^{\,*}}{\rm
dist}(\alpha^{\,*}_m, \beta^{\,*}_m)}\right)\geqslant$$
\begin{equation}\label{eq7C_1}
\geqslant \widetilde{c_n}\cdot \log\left
(1+\frac{\delta_2^{1/\alpha^{\,*}}}{(\widetilde{C})^{1/\alpha^{\,*}}|\varphi^{\,*}(x_m)-\varphi^{\,*}(y_m)|}\right)\,,
\end{equation}
де $\widetilde{c_n}>0$ -- деяка стала, яка залежить тільки від $n$ і
сталої $A^{\,*}_0$ з означення $QED$-області.

\medskip
Встановимо тепер верхню оцінку для $M(\Gamma_m).$ Зауважимо, що
відображення $F$ у~(\ref{eq16}) задовольняє
співвідношення~(\ref{eq2*A}) з функцією $\widetilde{Q(x)}=K_0\cdot
Q(x)$ замість $Q,$ де $K_0\geqslant 1$ -- стала квазіконформності
відображення $(\varphi^{\,*})^{\,-1}.$ Покладемо
$$\rho_m(y)= \left\{
\begin{array}{rr}
\frac{C^{\frac{1}{\alpha}}_1\left(C^{\,*}_1\right)^{\frac{1}{\alpha^2}}}{C^{\frac{1}{\alpha}}_0}\cdot
|f_m(x_m)-f_m(y_m)|^{-\frac{1}{\alpha^2}}, & y\in D^{\,\prime},\\
0,  &  y\not\in D^{\,\prime}\,.
\end{array}
\right. $$
Зауважимо, що $\rho_m$ задовольняє співвідношення~(\ref{eq1.4}) для
сім'ї кривих $f_m(\Gamma_m)$ в силу співвідношення~(\ref{eq22}).
Тоді за означення сім'ї ${\frak S}_{\delta, A, Q }$ ми отримаємо, що
\begin{equation}\label{eq34}
M(\Gamma_m)\leqslant
\frac{K_0\left(\frac{C^{\frac{1}{\alpha}}_1\left(C^{\,*}_1\right)^{\frac{1}{\alpha^2}}}{C^{\frac{1}{\alpha}}_0}\right)^n}
{|f_m(x_m)-f_m(y_m)|^{\frac{n}{\alpha^2}}} \cdot
\int\limits_{D^{\,\prime}} Q(y)\,dm(y)\,.
\end{equation}
З~(\ref{eq7C_1}) і (\ref{eq34}) випливає, що
$$\widetilde{c_n}\cdot \log\left
(1+\frac{\delta_2^{1/\alpha^{\,*}}}{(\widetilde{C})^{1/\alpha^{\,*}}|\varphi^{\,*}(x_m)-\varphi^{\,*}(y_m)|}\right)\leqslant
\frac{K_0\left(\frac{C^{\frac{1}{\alpha}}_1\left(C^{\,*}_1\right)^{\frac{1}{\alpha^2}}}{C^{\frac{1}{\alpha}}_0}\right)^n}
{|f_m(x_m)-f_m(y_m)|^{\frac{n}{\alpha^2}}}\cdot
\int\limits_{D^{\,\prime}} Q(y)\,dm(y)\,.$$
З останнього співвідношення, з огляду на гельдеревість відображення
$\varphi^{\,*}$ випливає, що
$$
|f_m(x_m)-f_m(y_m)|^{\frac{n}{\alpha^2}}\leqslant$$
\begin{equation}\label{eq33}
\leqslant K_0\Vert
Q\Vert_1\left(\frac{C^{\frac{1}{\alpha}}_1\left(C^{\,*}_1\right)^{\frac{1}{\alpha^2}}}{C^{\frac{1}{\alpha}}_0}\right)^n
\widetilde{c_n}\cdot \frac{1}{\log\left
(1+\frac{\delta_2^{1/\alpha^{\,*}}}{(\widetilde{C})
^{1/\alpha^{\,*}}|\varphi^{\,*}(x_m)-\varphi^{\,*}(y_m)|}\right)}\leqslant
\end{equation}
$$\leqslant \Vert
Q\Vert_1 M_0\cdot  \frac{\Vert
Q\Vert_1}{\log\left(1+\frac{\delta}{2|x_m-y_m|^{\alpha^{\,*}}}\right)}\,,$$
для деякої сталої $M_0>0,$ бо за правилом Лопіталя
$\log\left(1+\frac{1}{nt}\right)\sim\log\left(1+\frac{1}{kt}\right)$
при $t\rightarrow+0$ для різних фіксованих $k, n> 0.$ Оскільки при
кожному натуральному $m\in {\Bbb N}$ число $|f_m(x_m)-f_m(y_m)|$ є
меншим за одиницю, то
$|f_m(x_m)-f_m(y_m)|^{\frac{n}{\alpha^2}}\geqslant
|f_m(x_m)-f_m(y_m)|^{m}$ для достатньо великих $m.$ Тоді
з~(\ref{eq33}) випливає, що
$$
|f_m(x_m)-f_m(y_m)|^m\leqslant
|f_m(x_m)-f_m(y_m)|^{\frac{n}{\alpha^2}}\leqslant \Vert Q\Vert_1
M_0\cdot  \frac{\Vert
Q\Vert_1}{\log\left(1+\frac{\delta}{2|x_m-y_m|^{\,\alpha^{\,*}}}\right)}
\leqslant$$$$\leqslant
\Vert Q\Vert_1 m\cdot  \frac{\Vert
Q\Vert_1}{\log\left(1+\frac{\delta}{2|x_m-y_m|^{\,\alpha^{\,*}}}\right)}\,.$$
Останнє співвідношення суперечить припущенню~(\ref{eq2E}), що і
завершує доведення.~$\Box$
\end{proof}

\medskip
\begin{remark}\label{rem1}
Теорема~\ref{th3} залишається справедливою, якщо в ній одна, або дві
області $D$ і $D^{\,\prime}$ є, відповідно, локально
$\alpha^{\,*}$-гельдерово та $\alpha$-гельдерово еквівалентними до
півкулі на своїй межі. В цьому випадку будемо мати
співвідношення~(\ref{eq28}), у якому показник $\alpha^{\,*}$
відповідає <<гельдеровій еквівалентності>> $D,$ а показник $\alpha$
-- області $D^{\,\prime}.$ Доведення цього твердження з невеликими
відмінностями повторює доведення теореми~\ref{th3}, і тому не
наводиться.
\end{remark}

\medskip
{\bf 6. Прості кінці.} Означення простого кінця, яке
використовується нижче, може бути знайдено в~\cite{KR$_1$}. Тут і
далі $\overline{D}_P$ позначає поповнення області $D$ її простими
кінцями, а $E_D=\overline{D}_P\setminus D$ -- множина всіх простих
кінців в $D.$ Говоримо, що обмежена область $D$ в ${\Bbb R}^n$ {\it
регулярна}, якщо $D$ може бути квазіконформно відображена на область
з локально квазіконформною межею, замикання якої є компактом в
${\Bbb R}^n,$ крім того, кожен простий кінець $P\subset E_D$ є
регулярним. Зауважимо, що замикання $\overline{D}_P$ регулярної
області $D$ є {\it метризовним}, при цьому, якщо $g:D_0\rightarrow
D$ -- квазіконформне відображення області $D_0$ з локально
квазіконформною межею на область $D,$ то для $x, y\in
\overline{D}_P$ покладемо:
\begin{equation}\label{eq1A}
\rho(x, y):=|g^{\,-1}(x)-g^{\,-1}(y)|\,,
\end{equation}
де для $x\in E_D$ елемент $g^{\,-1}(x)$ розуміється як деяка (єдина)
точка межі $D_0,$ коректно визначена з огляду
на~\cite[теорема~4.1]{Na}. Має місце наступний результат.

\medskip
\begin{theorem}\label{th4}
{\sl\,Нехай $Q\in L^1(D^{\,\prime}),$ а $D$ є регулярною областю.

\medskip
\textbf{1)} Якщо $D^{\,\prime}$ є областю з локально квазіконформною
межею, то будь-яке відображення $f\in {\frak S}_{\delta, A, Q}(D,
D^{\,\prime})$ продовжується до відображення
$f:\overline{D}_P\rightarrow\overline{D^{\,\prime}},$ при цьому, для
кожної точки $P_0\in E_D$ знайдуться окіл $U$ цієї точки у
метричному просторі $(\overline{D}_P, \rho)$ і сталі $C_n=C(n, A, D,
D^{\,\prime})>0,$ $0<\alpha=\alpha(n, A, D, D^{\,\prime})\geqslant
1$ і $0<\alpha^{\,*}=\alpha^{\,*}(n, A, D, D^{\,\prime})\geqslant 1$
такі, що
\begin{equation}\label{eq2C_2}
|\overline{f}(P_1)-\overline{f}(P_2)|^{\frac{n}{\alpha^2}}\leqslant\frac{C_n\cdot
(\Vert
Q\Vert_1)^{1/n}}{\log^{1/n}\left(1+\frac{r_0}{\rho^{\alpha^{\,*}}(P_1,
P_2)}\right)}
\end{equation}
для всіх $P_1, P_2\in U,$ де $\Vert Q\Vert_1$ -- норма функції $Q$ в
$L^1(D).$

\medskip
\textbf{2)} Якщо $D^{\,\prime}$ є локально $\alpha$-гельдерово
еквівалентною до півкулі на своїй межі в сенсі
співвідношення~(\ref{eq25}), то будь-яке відображення $f\in {\frak
S}_{\delta, A, Q}(D, D^{\,\prime})$ продовжується до відображення
$f:\overline{D}_P\rightarrow\overline{D^{\,\prime}},$ при цьому, для
кожної точки $P_0\in E_D$ знайдуться окіл $U$ цієї точки у
метричному просторі $(\overline{D}_P, \rho)$ і сталі $C_n=C(n, A, D,
D^{\,\prime})>0,$ $r_0=r_0(n, A, x_0, D)>0$ такі, що виконано
співвідношення~(\ref{eq2C_2}).
 }
\end{theorem}

\medskip
\begin{proof} Нехай $f\in {\frak S}_{\delta, A, Q}(D, D^{\,\prime}).$
Достатньо обмежитися випадком $P_1, P_2\in U\cap D.$
Оскільки $D$ -- регулярна область, існує квазіконформне відображення
$g^{\,-1}$ області $D$ на область $D_0$ з локально квазіконформною
межею, причому, за означенням метрики $\rho$ в~(\ref{eq1A}),
\begin{equation}\label{eq18}
\rho(P_1, P_2):=|g^{\,-1}(P_1)-g^{\,-1}(P_2)|\,. \
\end{equation}
Розглянемо допоміжне відображення
\begin{equation}\label{eq19}
F(x)=(f\circ g)(x)\,,\quad x\in D_0\,.
\end{equation}
Оскільки відображення $g^{\,-1}$ є квазіконформним, то існує стала
$1\leqslant K_1<\infty$ така, що
\begin{equation}\label{eq17}
\frac{1}{K_1}\cdot M(\Gamma)\leqslant M(g^{\,-1}(\Gamma))\leqslant
K_1\cdot M(\Gamma)
\end{equation}
для будь-якої сім'ї кривих $\Gamma$ в $D_0.$ З огляду на
нерівності~(\ref{eq17}) та з урахуванням того, що $f$ задовольняє
співвідношення~(\ref{eq2*A}), ми отримаємо, що також $F$ задовольняє
співвідношення~(\ref{eq2*A}) з новою функцією $\widetilde{Q}(x):=
K_1\cdot Q(x).$ Крім того, оскільки $g$ -- фіксоване відображення,
яке є гомеоморфізмом, то $h(F^{\,-1}(A),
\partial D)\geqslant~\delta_0>0,$ де $\delta_0>0$ -- деяке
фіксоване число. Оскільки $g$ -- квазіконформне відображення, воно є
локально гельдеровим з деяким показником $\alpha^{\,*}.$

\medskip
Розглянемо випадок~\textbf{1)}. Тоді до відображення $F$ можна
застосувати теорему~\ref{th3}. Застосовуючи цю теорему, ми
отримаємо, що для будь-якої точки $x_0\in D_0$ знайдуться окіл
$U^{\,*}$ цієї точки і сталі $C^*_n=C(n, A, D_0, D^{\,\prime})>0,$
$0<\alpha^{\,*}=\alpha^{\,*}(n, A, D_0, D^{\,\prime})\geqslant 1$ і
$0<\alpha=\alpha(n, A, D_0, D^{\,\prime})\geqslant 1$ такі, що
\begin{equation}\label{eq2C_3}
|F(x)-F(y)|^{\frac{n}{\alpha^2}}\leqslant\frac{C^*_nK_1\cdot \Vert
Q\Vert_1}{\log\left(1+\frac{\delta_0}{2|x-y|^{\alpha^{\,*}}}\right)}
\end{equation}
для всіх $x, y\in U^{\,*}\cap D_0,$
де $\Vert Q\Vert_1$ -- норма функції $Q$ в $L^1(D).$ Нехай
$U:=g(U^{\,*}),$ $P_0:=g(x_0).$ Тоді за означенням $U$ є околом
простого кінця $P_0\in E_D.$ Якщо $P_1, P_2\in  D_P\cap U ,$ то
$P_1=g(x)$ і $P_2=g(y)$ для деяких $x, y \in U^{\,*}\cap D_0.$ З
огляду на співвідношення~(\ref{eq2C_3}) та враховуючи те, що
$|x-y|=|g^{\,-1}(P_1)-g^{\,-1}(P_2)|=\rho(P_1, P_2),$ ми отримаємо,
що
$$|F(g^{\,-1}(P_1))-F(g^{\,-1}(P_2))|^{\frac{n}{\alpha^2}}\leqslant
\frac{C^*_nK_1\cdot \Vert
Q\Vert_1}{\log\left(1+\frac{\delta_0}{2\rho^{\alpha^{\,*}}(P_1,
P_2)}\right)}\,,$$
або, з огляду на~(\ref{eq19}),
$$|f(P_1)-f(P_2)|^{\frac{n}{\alpha^2}}
\leqslant\frac{C^*_nK_1\cdot \Vert
Q\Vert_1}{\log\left(1+\frac{\delta_0}{2\rho^{\alpha^{\,*}}(P_1,
P_2)}\right)}\,.$$
Останнє співвідношення і є бажаним.

\medskip
Щодо випадку~\textbf{2)}, із зауваження~\ref{rem1} випливає
співвідношення~(\ref{eq2C_3}) з конкретними $\alpha^{\,*}$ і
$\alpha.$ Решта доведення залишиться незмінним.~$\Box$
\end{proof}


КОНТАКТНА ІНФОРМАЦІЯ

\medskip
\noindent{{\bf Марія Вікторівна Андрощук} \\
Житомирський державний університет ім.\ І.~Франко\\
вул. Велика Бердичівська, 40 \\
м.~Житомир, Україна, 10 008 \\
e-mail: mariaandroschuk28082003@gmail.com }

\medskip
\noindent{{\bf Олександр Петрович Довгопятий} \\
Житомирський державний університет ім.\ І.~Франко\\
вул. Велика Бердичівська, 40 \\
м.~Житомир, Україна, 10 008 \\
e-mail: alexdov1111111@gmail.com}

\medskip
\noindent{{\bf Наталія Сергіївна Ількевич} \\
Житомирський державний університет ім.\ І.~Франко\\
вул. Велика Бердичівська, 40 \\
м.~Житомир, Україна, 10 008 \\
e-mail: ilkevych1980@gmail.com}

\medskip
\noindent{{\bf Євген Олександрович Севостьянов} \\
{\bf 1.} Житомирський державний університет ім.\ І.~Франко\\
кафедра математичного аналізу, вул. Велика Бердичівська, 40 \\
м.~Житомир, Україна, 10 008 \\
{\bf 2.} Інститут прикладної математики і механіки
НАН України, \\
вул.~генерала Батюка, 19 \\
м.~Слов'янськ, Україна, 84 116\\
e-mail: esevostyanov2009@gmail.com}

\end{document}